 \documentclass[11pt,a4paper]{amsart}
\usepackage{amssymb,amsmath,amsthm}
\usepackage[left=2.7cm, right=2.7cm, marginparwidth=2.2cm, textheight = 24cm ]{geometry}
\usepackage{fullpage}

\usepackage[colorlinks=true,linkcolor=blue, citecolor=red]{hyperref}
\usepackage{eucal, lmodern}
\usepackage{color}
\usepackage{stackrel}
\usepackage{graphicx}
\newcommand{\ie}{\emph{i.e.}}
\newcommand{\eg}{\emph{e.g.}}

\newcommand{\etal}{\emph{et al.}}

\newcommand{\Real}{\mathbb{R}}

\newcommand{\dom}{\mathop{\mathrm{dom}}\nolimits}

\newcommand{\eps}{\varepsilon}
\newcommand{\sii}{L^2}

\newcommand{\der}{\mathrm{d}}

\newtheorem{Theorem}{Theorem}

\theoremstyle{definition}

\def\OMIT#1{}
\usepackage[normalem]{ulem}
\definecolor{DarkGreen}{rgb}{0,0.5,0.1} 
\definecolor{DarkBlue}{rgb}{0,0.1,0.5}

\newcommand\soutD{\bgroup\markoverwith
{\textcolor{DarkGreen}{\rule[.5ex]{2pt}{1pt}}}\ULon}
\newcommand\soutP{\bgroup\markoverwith
{\textcolor{blue}{\rule[.5ex]{2pt}{1pt}}}\ULon}
\newcommand{\Hm}[1]{\leavevmode{\marginpar{\tiny%
$\hbox to 0mm{\hspace*{-0.5mm}$\leftarrow$\hss}%
\vcenter{\vrule depth 0.1mm height 0.1mm width \the\marginparwidth}%
\hbox to
0mm{\hss$\rightarrow$\hspace*{-0.5mm}}$\\
\relax\raggedright #1}}}

\newcommand\Omg{\Omega}
\newcounter{counter_a}
\newenvironment{myenum}{\begin{list}{{\rm(\roman{counter_a})}}%
{\usecounter{counter_a}
\setlength{\itemsep}{1.ex}\setlength{\topsep}{0.8ex}
\setlength{\leftmargin}{5ex}
\setlength{\labelwidth}{5ex}}}{\end{list}}
\def\dR{{\mathbb R}}
\newcommand\frh{\mathfrak h}
\newcommand\nb{\nabla}

\def\s{\sigma}
\def\lm{\lambda}

\def\dN{{\mathbb N}}
\def\arr{\rightarrow}

\newtheorem{thm}{Theorem}[section]
\newtheorem{dfn}[thm]{Definition}
\newtheorem{prop}[thm]{Proposition}
\newtheorem{rem}[thm]{Remark}
\def\G{\Gamma}
\newcommand\ov{\overline}
\def\p{\partial}
\def\sfK{{\mathsf K}} 
\def\sfB{{\mathsf B}}
\def\omg{\omega}

\newcommand\dl{\delta}
\def\aa{\alpha}
\newcommand\sm{\setminus} 
\def\sfT{{\mathsf T}}
\newcommand\sfS{{\mathsf S}}

\title{Quasi-conical domains with embedded eigenvalues}
\author[D.~Krej\v{c}i\v{r}\'{i}k]{David Krej\v{c}i\v{r}\'{i}k}
\address{(D.~Krej\v{c}i\v{r}\'{i}k) Department of Mathematics\\ Faculty of Nuclear Sciences and Physical
	Engineering\\
	Czech Technical University in Prague\\ Trojanova 13, 120 00, Prague, Czech
	Republic \newline 
	{E-mail: david.krejcirik@fjfi.cvut.cz}}
\author[V.~Lotoreichik]{Vladimir Lotoreichik}
\address{(V.~Lotoreichik)
	Department of Theoretical Physics\\
	Nuclear Physics Institute, Czech Academy of Sciences, 
	25068 \v{R}e\v{z}, Czech Republic \newline
	E-mail: {lotoreichik@ujf.cas.cz}
}


\begin{document}
%
%

\begin{abstract}
	The spectrum of the Dirichlet Laplacian on any quasi-conical open set coincides with the non-negative semi-axis.
	We show that there is a connected quasi-conical open set such that the
	respective Dirichlet Laplacian has a
	positive (embedded) eigenvalue. 
	This open set is constructed as the tower of cubes of growing size connected by  windows of vanishing size. Moreover, we show that the sizes of the windows in this construction can be chosen so  that the absolutely continuous spectrum of the Dirichlet Laplacian is empty.
\end{abstract}

\maketitle
\section{Introduction}
An open set $\Omega \subset \Real^d$ with any $d \geq 2$
is said to be \emph{quasi-conical}
if it contains an arbitrarily large ball.
We consider the self-adjoint Dirichlet Laplacian $-\Delta_{\rm D}^\Omg$ on $\Omg$ acting in the Hilbert space $L^2(\Omg)$. 
It is easy to see that 
\begin{equation}\label{spectrum}
  \sigma(-\Delta_{\rm D}^\Omega) = [0,\infty)
\end{equation}
for any quasi-conical open set.

Sufficient conditions to ensure an absence of eigenvalues 
of the Dirichlet Laplacian on quasi-conical sets
are available in the literature. In particular, Rellich proved in~\cite{Rellich_1943} that there are no embedded eigenvalues for the Dirichlet Laplacian on the exterior of a compact set under
the assumption that this exterior is connected.
Jones~\cite{Jones_1953} excluded embedded eigenvalues for certain other classes of quasi-conical open sets such as domains conical at infinity.
More recently, D'Ancona and Racke ~\cite{D'Ancona-Racke_2012} excluded embedded eigenvalues
by imposing a repulsive-type condition on the geometry of the boundary
of tubular-type quasi-conical sets.
Finally, Bonnet-Ben Dhia \etal\
\cite{Bonnet-Fliss-Hazard-Tonnoir_2016}
excluded the existence of eigenvalues in non-convex conical sectors.
(The critical case of the Robin Laplacian in a half-space 
is considered in~\cite{CK2}.)

\begin{figure}[h!]
\begin{center}
\includegraphics[width=0.9\textwidth]{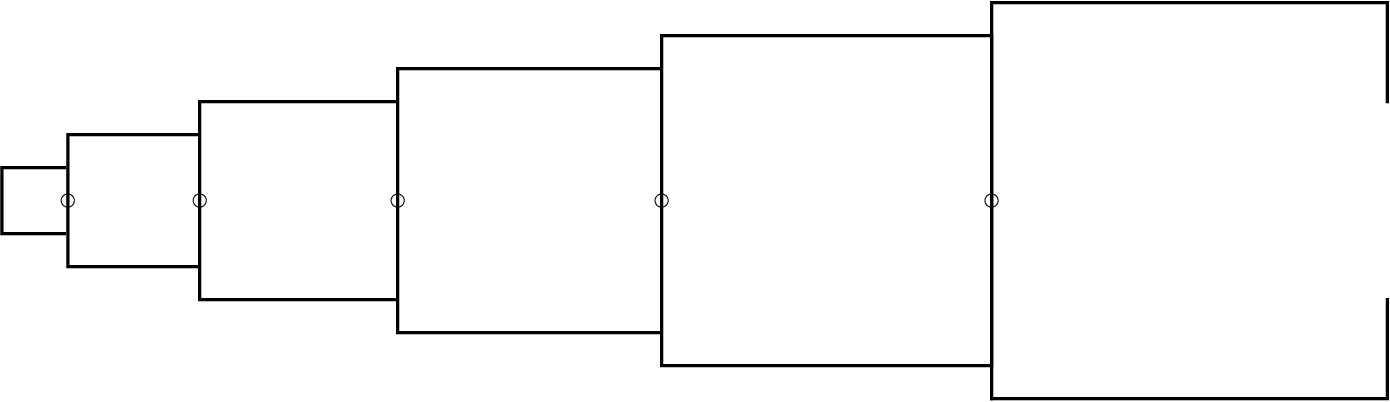}
\end{center}
\label{fig:tower}
\caption{The disconnected quasi-conical open set $\Omg$ with places for the small holes to be dug indicated by circles.}
\end{figure}

On the other hand,
it is easy to construct examples of disconnected quasi-conical open sets 
for which there are (embedded) eigenvalues in $[0,\infty)$
(\eg, take the union of any bounded set and its exterior).
The question whether there are examples of quasi-conical \emph{domains}
(\ie~\emph{connected} open sets)
for which there are eigenvalues seems not to have
been considered in the literature. 
Our goal is to fill in this gap.
The main result of the present paper reads as follows.
\begin{Theorem}\label{thm:main}
\ 
\begin{myenum}	
\item Given any positive number~$\lambda$, 
there exists a connected quasi-conical open set such that 
\begin{equation*}
  \lambda \in \sigma_\mathrm{p}(-\Delta_{\rm D}^{\Omega})
  \,.
\end{equation*}
\item The connected quasi-conical open set in {\rm (i)} can be constructed so that 
\begin{equation*}
	\sigma_\mathrm{ac}(-\Delta_{\rm D}^{\Omg}) = \varnothing\,.
\end{equation*}
\end{myenum}
%
\end{Theorem}
Our proof relies on an explicit construction, 
which we believe is of independent interest. The desired connected quasi-conical open set is constructed as the tower of cubes of growing size connected by windows of vanishing size (see Figure~\ref{fig:tower}). 
We show that the sizes of the windows connecting the cubes can be chosen so small that we get an embedded eigenvalue and that the absolutely continuous spectrum turns out to be empty. Our result is partially of qualitative nature, because we have no explicit control on the sizes of the windows, which ensure the desired property.

In order to show the existence of an embedded eigenvalue, 
we start from the tower made of disjoint cubes. The respective Dirichlet Laplacian clearly has an embedded eigenvalue. Then we
start to subsequently open the windows
choosing their sizes so small that this embedded eigenvalue and the respective eigenfunction do not change much (in a proper sense). In the course of this construction we get a family of eigenfunctions weakly convergent to a function, which turns out to be an eigenfunction  corresponding to an embedded eigenvalue for the tower of cubes with all the windows opened. In this argument we rely on perturbation results for Dirichlet eigenvalues and eigenfunctions under variation of the open set~\cite{Daners}.

The absence of the absolutely continuous spectrum is shown via a perturbation type argument.
Clearly, for the disjoint cubes the absolutely continuous
spectrum of the Dirichlet Laplacian is empty.
Adapting the technique of~\cite{HSS91},
we show that provided the sizes of the windows are sufficiently small, a suitable resolvent power difference between the Dirichlet Laplacian on a tower of cubes with opened windows and the Dirichlet Laplacian on the disjoint cubes belongs to the trace class.  The claim follows from the stability of the absolutely continuous spectrum under such a perturbation.

In the classification of Glazman~\cite{Glazman} (see also~\cite[Section X.6.1]{Edmunds-Evans}) there are three types of open sets: 
\emph{quasi-conical} (\ie\ presently considered), 
\emph{quasi-cylindrical} (\ie\ not quasi-conical 
but containing infinitely many pairwise disjoint balls of the same radius), 
and \emph{quasi-bounded} (\ie\ neither quasi-conical nor quasi-cylindrical). 
Quasi-cylindrical open sets with embedded eigenvalues are constructed in a large number of papers. 
Behind these constructions there is typically a symmetry based argument.
For example, it follows from the considerations in
\cite{ESS1,ABGM} (see also~\cite{DLR12}) 
that there is an embedded eigenvalue in the essential spectrum
of the Dirichlet Laplacian in an infinite Swiss-cross waveguide.
Other examples are given in \cite{Witsch_1990a,Witsch_1990b}.
The construction in the present paper for quasi-conical domains 
is significantly different and not based on symmetry. 
The Dirichlet Laplacian on a quasi-bounded open set typically has purely discrete spectrum, but not necessarily if the set is not bounded.%
\footnote{
In passing, our present construction can be considered
as a dual geometry to the celebrated domain ``rooms and passages''
due to Fraenkel~\cite{Fraenkel_1979},
which exhibits an essential spectrum for the Neumann Laplacian on a bounded domain.}

To the best of our knowledge,
as in the present situation of quasi-conical domains, 
embedded eigenvalues for unbounded quasi-bounded open sets 
are not studied in the literature.

The structure of the paper is as follows.
First, in Section~\ref{sec:pre} we provide preliminary material that is used throughout the paper. In Subsection~\ref{ssec:DirLap} we rigorously introduce the Dirichlet Laplacian on a quasi-conical open set 
and evoke its spectral property~\eqref{spectrum}.
Next, in Subsection~\ref{ssec:Daners} we recall general results on perturbation of
eigenvalues, eigenfunctions, and the resolvent of the Dirichlet Laplacian under variation of the open set. Further, in Subsection~\ref{ssec:sk} we recall the definitions of the singular values for a compact operator and of the trace class.
In Section~\ref{sec:proof1} we prove Theorem~\ref{thm:main}\,(i) on the existence of embedded eigenvalues. Theorem~\ref{thm:main}\,(ii) on the absence of the absolutely continuous spectrum is proved in Section~\ref{sec:proof2}. The paper is complemented by Section~\ref{sec:discussion} with a discussion of a related open question. 

\section{Preliminaries}\label{sec:pre}
\subsection{Dirichlet Laplacian on a (quasi-conical) open set}\label{ssec:DirLap}
Let $\Omg\subset\dR^d$ ($d\ge 2$) be an open set. 
The self-adjoint Dirichlet Laplacian $-\Delta_{\rm D}^\Omg$ in $L^2(\Omg)$ can be rigorously introduced using the first representation theorem~\cite[Thm.~10.7]{S12} via the quadratic form
\[
	\frh_{\rm D}^\Omg[u] := \int_\Omg|\nb u|^2\der x\,,
	\qquad\dom\frh_{\rm D}^\Omg := H^1_0(\Omg)\,.
\]
%
The operator $-\Delta_{\rm D}^\Omg$ can be characterised as
\[
	-\Delta_{\rm D}^{\Omg}u = -\Delta u\,,\qquad \dom(-\Delta_{\rm D}^{\Omg}) = \big\{u\in H^1_0(\Omg)\colon \Delta u \in L^2(\Omg)\big\}\,.
\]

Since the Dirichlet Laplacian is a non-negative operator, 
we have $\s(-\Delta_{\rm D}^\Omg) \subseteq [0,\infty)$
for any open set~$\Omega$.
Now, assume additionally that $\Omg$ is a quasi-conical open set
(recall that an open set is called quasi-conical if it contains a ball of arbitrarily large radius). 
Then, using quasi-modes based on plane waves modulated 
by suitable cut-off functions in the enlarging balls 
(see, \eg, \cite[Thm.~X.6.5]{Edmunds-Evans}),
the opposite inclusion of~\eqref{spectrum} follows.
In particular, the spectrum of the Dirichlet Laplacian 
on any quasi-conical open set is purely essential.

\subsection{Dirichlet Laplacian on varying domains}\label{ssec:Daners}
In this subsection we recall some convergence results for the Dirichlet Laplacian on varying domains contained in~\cite{Daners}
(see also \cite{Daners_2008, RT,Stollmann}). 
We begin by introducing a type of convergence for a sequence of open sets.
\begin{dfn}[{\cite[Def.~1.1]{Daners}}]
	\label{dfn:Mosco}
	If $\Omg_n,\Omg\subset\dR^d$ are such that
	\begin{myenum}
		\item the weak limit points of every sequence $u_n\in H^1_0(\Omg_n)$, $n\in\dN$, in $H^1(\dR^d)$ are in $H^1_0(\Omg)$,
		\item for every $\varphi\in H^1_0(\Omg)$ there exist $\varphi_n\in H^1_0(\Omg_n)$ such that $\varphi_n\arr\varphi$ in $H^1(\dR^d)$,
	\end{myenum}	
	then we write $\Omg_n\arr\Omg$ and say that $\Omg_n$ converges to $\Omg$ 
	\emph{in the sense of Mosco}.
\end{dfn}
We remark that in Definition~\ref{dfn:Mosco} the domains $\Omg_n,\Omg$ are not assumed to be always connected.
In general, verifying the conditions of Definition~\ref{dfn:Mosco} can be a difficult problem. However, in the case that the sequence $\Omg_n$ approximates $\Omg$ monotonously from inside or outside the convergence in the sense of Mosco can be easily checked. In the present paper we only use the approximation from outside.
(For any  set $S\subset\dR^d$ we denote by ${\rm int}(S)$ its interior.)

\begin{prop}[{\cite[Prop.~7.4 and Thm.~7.5]{Daners}, \cite[Sec.~5.4] {Daners_2008}\label{prop:Mosco}}]
	Suppose that $\Omg_n\supset\Omg_{n+1}\supset\Omg$ for all $n\in\dN$, and that ${\rm int}\,\big(\cap_{n\in\dN}\Omg_n\big) = \Omg$. If $H^1_0(\Omg) = H^1_0(\Omg\cup\G)$, where
	\[
		\G := \bigcap_{n\in\dN}\left(
		\ov{\bigcup_{k\ge n}(\Omg_k\cap\p\Omg)}\right)\subset\p\Omg\,,
	\]	%
	then $\Omg_n\arr\Omg$.
\end{prop}
The convergence of the open sets in the sense of Mosco is convenient in the study of convergence of the respective Dirichlet Laplacians. In order to state this convergence result, we define for an open set $\Omg\subset\dR^d$ the restriction mapping $r_\Omg$ and the mapping $i_\Omg$, which extends the function from $\Omg$ by zero 
to the whole $\dR^d$:
\[
\begin{aligned}
&r_\Omg\colon L^2(\dR^d)\arr L^2(\Omg), \qquad &r_\Omg u &:= u|_{\Omg},\\
&i_\Omg\colon L^2(\Omg)\arr L^2(\dR^d), \qquad &i_\Omg u& := u\oplus 0.
\end{aligned}
\]
\begin{prop}[{\cite[Cor.~4.7]{Daners}}]
	\label{prop:norm_resolvent}
	Suppose that $\Omg_n,\Omg$ are contained in a fixed bounded set and that $\Omg_n\arr\Omg$. Then
	\[
		\left\|i_{\Omg_n}(-\Delta_{\rm D}^{\Omg_n}+1)^{-1}r_{\Omg_n} - 
		i_\Omg(-\Delta_{\rm D}^{\Omg}+1)^{-1}r_\Omg\right\|\arr0\,,\qquad n\arr\infty\,.
	\]
\end{prop}
In the present paper, we apply the above proposition only under the additional assumption that $\ov{\Omg_n}$ coincides with $\ov\Omg$ for all $n\in\dN$. In this special case, 
$-\Delta^{\Omg_n}_{\rm D}$ and $-\Delta^\Omg_{\rm D}$ act in the same Hilbert space $L^2(\Omg_n) = L^2(\Omg)$ and	
 Proposition~\ref{prop:norm_resolvent} implies that
$-\Delta_{\rm D}^{\Omg_n}$ converges in the norm-resolvent sense to $-\Delta_{\rm D}^\Omg$ as $n\arr\infty$.  
This norm-resolvent convergence particularly implies 
the convergence of the respective spectral projections.

\subsection{Singular values and the trace class}\label{ssec:sk}
In this subsection we recall the concepts of singular values and the trace class. Further details can be found in the monographs~\cite{Gohberg-Krein_1969, Simon05}.

\begin{dfn}\label{dfn:sk}
For a compact operator $\sfK$
in a Hilbert space the eigenvalues  $\{s_k(\sfK)\}_{k\ge1}$ of the self-adjoint non-negative compact operator $(\sfK^*\sfK)^{1/2}$ enumerated in the non-increasing order and repeated with multiplicities taken into account are called the \emph{singular values} of $\sfK$.
\end{dfn}
For a compact operator $\sfK$ and a bounded closed operator $\sfB$ acting in a Hilbert space the products $\sfB\sfK$ and $\sfK\sfB$ are compact and
by~\cite[Thm.~1.6]{Simon05}
\begin{equation}\label{eq:sk1}
s_k(\sfB\sfK) \le \|\sfB\| \, s_k(\sfK)\quad\text{and}\quad
s_k(\sfK\sfB) \le \|\sfB\| \, s_k(\sfK).
\end{equation}
%
\begin{dfn}\label{dfn:traceclass}
A compact operator $\sfK$ is said to belong to 
the \emph{trace class} 
if $\displaystyle \sum_{k=1}^\infty s_k(\sfK) < \infty$. 
\end{dfn}
According to~\cite[Chap.~2]{Simon05}, the trace class is a Banach space with the norm given by
\begin{equation}\label{eq:tracenorm}
\|\sfK\|_1 := \sum_{k=1}^\infty s_k(\sfK).
\end{equation}

\section{Embedded eigenvalues}\label{sec:proof1}
In this section we prove 
part~(i) of
Theorem~\ref{thm:main} on the existence of a connected quasi-conical open set with an embedded eigenvalue. The argument is performed via an explicit construction.
For clarity, we divide the proof into several steps.

\subsubsection*{Disconnected set}
Let~$Q_a(x)\subset\dR^d$ denote an open cube with sides of half-width $a>0$
centred at a point $x \in \Real^d$.
With a slight abuse of notation we will occasionally interpret $t\in\dR$ as $(t,0,0,\dots,0)\in\dR^d$. We aim at constructing a tower of cubes of growing size, in which each next cube is placed ``on the top" of the previous one. To this aim, let us introduce the sequence of centres $0 < x_1 < x_2 < \dots < x_n <\dots$ such that
the associated sequence of half-widths $a_1 :=x_1$, $a_n :=x_n - x_{n-1}-a_{n-1}$ ($n\ge 2$) satisfies $a_n > 0$ for all $n\in\dN$ and $a_n \arr\infty$ as $n\arr\infty$. At this stage we fix only $x_1> 0$ while the sequence $\{x_n\}_{n\ge2}$ is not fixed and we will specify it in the course of the construction.
Let $\lm_1$ be the lowest (simple) eigenvalue of $-\Delta_{\rm D}^{Q_{a_1}(x_1)}$. 
Let~$\psi_1$ denote a corresponding eigenfunction,
normalised to~$1$ in $\sii(Q_{a_1}(x_1))$. 

Define a (disconnected) quasi-conical open set~$\Omega$
as a disjoint union of enlarging cubes:
\begin{equation}
  \Omega := \bigcup_{n=1}^\infty \omega_n 
  \qquad \mbox{with} \qquad
  \omega_n := Q_{a_n}(x_n)
  \,.
\end{equation}
Then (extending~$\psi_1$ by zero to~$\Omega$),
it is easy to see that~$\psi_1$ is an eigenfunction of $-\Delta_{\rm D}^\Omega$
corresponding to~$\lambda_1$. 
Thus, we have constructed a disconnected quasi-conical open set with an embedded eigenvalue.

To construct a connected open set based on~$\Omega$,
we dig small holes at the places where the disjoint cubes~$\omega_n$ touch (see Figure~\ref{fig:tower}).
A careful attention is needed in order to ensure that this construction
does not eliminate the embedded eigenvalue. 

\subsubsection*{1st interconnection}
Let us set $\Omg_1 :=\omg_1$ and let us choose $x_2 > x_1$ so that $a_2 = x_2-x_1 -a_1 \ge a_1 +1$ and that $\lm_1$ is a simple eigenvalue of $-\Delta_{\rm D}^{\Omg_1\cup\omega_2}$.
For any positive $\delta_1 < a_1$, let us connect 
the first two cubes $\Omega_1 = \omega_1$ and~$\omega_2$ by defining
\begin{equation}
  \Omega_2 := \Omega_1 \cup \omega_2
  \cup Q_{\delta_1}(x_1+a_1)
  \,.
\end{equation}
It follows from Proposition~\ref{prop:Mosco} (with $\G = \{x_1+a_1\}$) that $\Omg_2$ converges to $\Omg_1\cup\omg_2$ as $\dl_1\arr0$ in the sense of Definition~\ref{dfn:Mosco}.
Hence, Proposition~\ref{prop:norm_resolvent} implies
that $-\Delta_{\rm D}^{\Omg_2}$ converges in the norm-resolvent sense to $-\Delta_{\rm D}^{\Omg_1\cup\omg_2}$ as $\dl_1\arr 0$ and combined with~\cite[Sec.~IV.3.5]{Kato} yields that 
given any (small) positive number~$\eps_1$, there exist a positive half-width~$\delta_1$, 
a simple eigenvalue~$\lambda_2$ of $-\Delta_{\rm D}^{\Omega_2}$ 
and a corresponding eigenfunction~$\psi_2$,
normalised to~$1$ in $\sii(\Omega_2)$,
such that  
\begin{equation}\label{p1}
  |\lambda_2-\lambda_1| \leq \eps_1 
  \qquad\mbox{and}\qquad
  \|\psi_2-\psi_1\|_{\sii(\Omega_2)} \leq \eps_1 
  \,.
\end{equation}
Using additionally the monotonicity of Dirichlet eigenvalues with respect to inclusion of open sets,
the first inequality in~\eqref{p1} implies
that we can choose $\dl_1$ so small that
\begin{equation}\label{l1}
  \lambda_1 - \eps_1 \leq \lambda_2 \leq \lambda_1 
  \,.
\end{equation}
Using the normalisations of~$\psi_1$ and~$\psi_2$,
together with the fact that $\psi_1=0$ in~$\omega_2$,  
the second inequality in~\eqref{p1} yields
\begin{equation}\label{r1}
  \|\psi_2\|_{\sii(\Omega_2\setminus\Omega_1)} \leq \eps_1
  \qquad\mbox{and}\qquad
  \|\psi_2\|_{\sii(\Omega_1)} \geq 1-\eps_1
  \,.
\end{equation}
\subsubsection*{2nd interconnection}
Similarly, let us choose $x_3 > x_2$
so that $a_3 = x_3-x_2-a_2 \ge a_2+1$ and that~$\lm_2$ 
(constructed in the previous step) is a simple eigenvalue of $-\Delta_{\rm D}^{\Omg_2\cup\omg_3}$.
For any $\delta_2 < a_2$, let us connect~$\Omega_2$ 
with the third cube~$\omega_3$ by defining
\begin{equation}
  \Omega_3 := \Omega_2 \cup \omega_3
  \cup Q_{\delta_2}(x_2+a_2)
  \,.
\end{equation}
Arguing as before, given any positive number~$\eps_2$, there exist a positive half-width~$\delta_2$,  
a simple eigenvalue~$\lambda_3$ of $-\Delta_{\rm D}^{\Omega_3}$, 
and a corresponding eigenfunction~$\psi_3$, 
normalised to~$1$ in $\sii(\Omega_3)$,
such that 
\begin{equation}\label{p2}
  |\lambda_3-\lambda_2| \leq \eps_2 
  \qquad\mbox{and}\qquad
  \|\psi_3-\psi_2\|_{\sii(\Omega_3)} \leq \eps_2
  \,.
\end{equation}
Using the monotonicity of Dirichlet eigenvalues and~\eqref{l1},
the first inequality in~\eqref{p2} implies
that we can choose $\dl_2$ so small that
\begin{equation}\label{l2}
  \lambda_1 - \eps_1 - \eps_2 \leq \lambda_3 \leq\lm_2\leq \lambda_1
  \,.
\end{equation}
From the second inequality in~\eqref{p2}
and the first inequality in~\eqref{r1}, 
we deduce
\begin{equation} 
  \|\psi_3\|_{\sii(\Omega_3\setminus\Omega_1)} 
  \leq \|\psi_3-\psi_2\|_{\sii(\Omega_3\setminus\Omega_1)} 
  + \|\psi_2\|_{\sii(\Omega_3\setminus\Omega_1)} 
  \leq \eps_2 + \eps_1
  \,.
\end{equation}
Consequently,
\begin{equation}
  1 = \|\psi_3\|_{\sii(\Omega_3)}
  \leq \|\psi_3\|_{\sii(\Omega_1)}
  + \|\psi_3\|_{\sii(\Omega_3\setminus\Omega_1)}
  \leq \|\psi_3\|_{\sii(\Omega_1)}
  + \eps_2 + \eps_1 
  \,.
\end{equation}
That is, 
\begin{equation}
  \|\psi_3\|_{\sii(\Omega_1)} \geq 1 - \eps_1 - \eps_2
  \,.
\end{equation}

\subsubsection*{$n$th interconnection}
Continuing by induction, 
given any natural number $n \geq 1$,
we choose $x_{n+1} > x_{n}$ so that $a_{n+1} = x_{n+1} -x_n -a_n \ge a_n+1$ and $\lm_n$ (constructed in the $(n-1)$st step) is a simple eigenvalue of $-\Delta_{\rm D}^{\Omg_n\cup\omg_{n+1}}$.
For any $\delta_{n} < a_n$,
let us connect~$\Omega_{n}$ 
with the $(n+1)$th cube $\omega_{n+1}$ by defining
\begin{equation}
  \Omega_{n+1} := \Omega_{n} \cup \omega_{n+1}
  \cup Q_{\delta_{n}}\big(x_n+a_n\big)
  \,.
\end{equation}
As before given any positive number~$\eps_n$, there exist a positive half-width~$\delta_n$, 
a simple eigenvalue~$\lambda_{n+1}$ of $-\Delta_{\rm D}^{\Omega_{n+1}}$ 
and a corresponding eigenfunction~$\psi_{n+1}$, 
normalised to~$1$ in $\sii(\Omega_{n+1})$,
such that  
\begin{equation}
  \lambda_1 - (\eps_1+\dots+\eps_n) \leq \lambda_{n+1} \leq\lm_n\leq \lambda_1\,,
\end{equation}
and that
\begin{equation}\label{eq:pn}  
  \|\psi_{n+1}-\psi_n\|_{L^2(\Omega_{n+1})} \le\eps_n\qquad\mbox{and}\qquad
  \|\psi_{n+1}\|_{\sii(\Omega_1)} \geq 1 -(\eps_1+\dots+\eps_n) 
  \,.  
\end{equation}
Given any positive number~$\eps$, 
and choosing $\eps_k := 2^{-k} \eps$ for every $k \geq 1$,
we get uniform estimates
\begin{equation}\label{uniform}
  \lambda_1 - \eps \leq \lambda_{n+1} \leq \lambda_1
  \qquad\mbox{and}\qquad
  \|\psi_{n+1}\|_{\sii(\Omega_1)} \geq 1 - \eps 
  \,.  
\end{equation}
Moreover, we get that the sequence $\{\psi_n\}_{n=1}^\infty$ is a Cauchy sequence in the Hilbert space $L^2(\Omega)$, where the eigenfunctions $\psi_n$ are viewed as functions in $L^2(\Omega)$ upon extension by zero. Indeed, for any $N\in\dN$ and any
integer $m,n \ge N$ with $m\ge n$, we get using the first inequality in~\eqref{eq:pn} and the triangle inequality for the norm that
\[
	\|\psi_m - \psi_n\|_{L^2(\Omg)} \le \sum_{k=n}^{m-1}
		\|\psi_{k+1} - \psi_k\|_{L^2(\Omg_{k+1})}
	\le \eps\sum_{k=n}^{m-1} 2^{-k} \le \frac{\eps}{2^{N-1}}\,.
\]
\subsubsection*{Connected set}
By using the procedure of connecting the individual cubes described above, 
we define a connected open set~$\Omega'$ by the limit
\begin{equation}
  \Omega' := \lim_{n \to \infty} \Omega_n 
  \,.  
\end{equation}
It remains to show that the embedded eigenvalue ``survives'' this limit.  
After the $n$th interconnection, we have the eigenvalue equation
\begin{equation}\label{weak}
  \forall \phi \in H_0^1(\Omega_{n+1})
  \,, \qquad
  \big(\nabla\phi,\nabla\psi_{n+1}\big)_{\sii(\Omega_{n+1})} 
  = \lambda_{n+1} \, \big(\phi,\psi_{n+1}\big)_{\sii(\Omega_{n+1})} 
  \,.  
\end{equation}
As usual, extending~$\psi_{n+1}$ and~$\phi$ by zero to the whole~$\Omega'$,
the inner product can be replaced by that of $\sii(\Omega')$. 
Choosing $\phi := \psi_{n+1}$ for the test function in~\eqref{weak}
and using the normalisation of~$\psi_{n+1}$,
we obtain  
\begin{equation} 
  \|\nabla \psi_{n+1}\|_{\sii(\Omega')}^2 
  = \lambda_{n+1} \leq \lambda_1
  \qquad\mbox{and}\qquad
  \|\psi_{n+1}\|_{\sii(\Omega')}^2 = 1 
  \,.  
\end{equation}
Here the uniform bound to~$\lambda_{n+1}$ follows from~\eqref{uniform}.
It follows that $\{\psi_{n+1}\}_{n=0}^\infty$ is a bounded sequence
in $H_0^1(\Omega')$.
Consequently, there exists a subsequence $\{\psi_{n_j+1}\}_{j=0}^\infty$
which is weakly converging to some $\psi_\infty$ in $H_0^1(\Omega')$.
Being a Cauchy sequence, $\{\psi_{n_j+1}\}_{j=0}^\infty$ is also strongly converging to $\psi_\infty$ in $L^2(\Omg')$.  
Inserting  this subsequence  into~\eqref{weak}, sending~$j$ to infinity, 
and taking into account that any $\phi\in C^\infty_0(\Omg')$ belongs to $H^1_0(\Omg_{n+1})$ for all $n\in\dN$ large enough, we obtain   
\begin{equation}\label{weak.bis}
  \forall \phi \in C^\infty_0(\Omg')
  \,, \qquad
  \big(\nabla\phi,\nabla\psi_{\infty}\big)_{\sii(\Omega')} 
  = \lambda_\infty \, \big(\phi,\psi_{\infty}\big)_{\sii(\Omega')} 
  \,,
\end{equation}
where the limit
\begin{equation}
  \lambda_\infty := \lim_{n\to\infty} \lambda_n 
\end{equation}
is well defined because $\{\lambda_{n}\}_{n=1}^\infty$
is monotonously decreasing by construction. 
Since $C^\infty_0(\Omg')$ is dense in $H^1_0(\Omg')$, we conclude that the equality in~\eqref{weak.bis} holds also for all $\phi \in H^1_0(\Omg')$.
Consequently, $\psi_\infty$ belongs to the domain of 
the operator $-\Delta_{\rm D}^{\Omega'}$
and solves the equation 
$-\Delta_{\rm D}^{\Omega'}\psi_\infty = \lambda_\infty\psi_\infty$
with $\lambda_1-\eps \leq  \lambda_\infty \leq \lambda_1$. 
By choosing $\eps < \lambda_1$
we get that $\lambda_\infty > 0$.
It remains to ensure that~$\psi_\infty$ 
is not identically equal to zero.
However, this follows from the 
fact that $\psi_\infty$ is the strong limit of a normalised sequence of functions. In particular, we infer that $\|\psi_\infty\|_{L^2(\Omega')} = 1$.

We remark that by rescaling the domain $\Omg'$ by a constant factor $\aa > 0$ we get the domain $\aa\Omg'$ and the respective Dirichlet Laplacian $-\Delta_{\rm D}^{\aa\Omg'}$ has an embedded eigenvalue $\aa^{-2}\lm_\infty$. Hence, we can construct a connected quasi-conical open set with an embedded eigenvalue at any prescribed position.  
\begin{rem}
It is essential for our approach to choose
$x_{n+1} > x_n$ in the $n$th step of the construction, so that $\lm_n$ is a simple eigenvalue in the spectrum $-\Delta_{\rm D}^{\Omg_n\cup\omg_{n+1}}$ in order to guarantee that properly normalized perturbed eigenfunction $\psi_{n+1}$ is close to $\psi_n$ in the sense of the first inequality in~\eqref{eq:pn}. Without preserving the simplicity of the eigenvalue $\lm_n$ upon adding the cube $\omg_{n+1}$, we can only claim that the spectral projections are close after digging a small hole in the $n$th interconnection, from which closeness of the eigenfunctions themselves in $L^2$-norm does not follow in general.
\end{rem}
\begin{rem}
	It is worth adding that as a consequence of the
	uniform estimate~\eqref{uniform} 
	on the norm of~$\psi_{n+1}$ in the $n$-independent bounded set~$\Omega_1$ with $\eps < 1$, we obtain the localisation estimate $\|\psi_\infty\|_{L^2(\Omg_1)} \ge 1-\eps$. By choosing $\eps > 0$ sufficiently small we can construct the connected quasi-conical open set $\Omega'$ such that the operator $-\Delta_{\rm D}^{\Omg'}$ has an embedded eigenvalue with the eigenfunction mostly localised in the first cube $\Omg_1$. 
\end{rem}
\begin{rem}
The same type of examples 
of quasi-conical domains with embedded eigenvalues
can be given for other cell domains than cubes.
Moreover, the interconnections obtained by digging the holes 
for flat boundaries can be replaced for general cell domains 
by using thin tubes as connections.
\end{rem}

%
\section{Absence of the absolutely continuous spectrum}\label{sec:proof2}

In this section we prove part~(ii) of Theorem~\ref{thm:main}
by showing that  
by choosing the sizes of the windows 
$\{\dl_n\}_{n\ge1}$ sufficiently small the absolutely continuous spectrum of the operator $-\Delta_{\rm D}^{\Omg'}$ becomes empty. In order to prove this result we adapt the technique developed in~\cite{HSS91}. 
As in Section~\ref{sec:proof1}, we subsequently open the windows between the cubes. 
The half-widths of the windows $\{\dl_n\}_{n\ge1}$ will be determined so as to satisfy simultaneously the conditions in Section~\ref{sec:proof1} and some new additional conditions. 

Let the sequence of the bounded domains $\{\Omg_n\}_{n\in\dN}$ be constructed as in Section~\ref{sec:proof1}. We intend to modify the sizes of the windows $\{\dl_n\}_{n\ge1}$ in order to satisfy some additional assumptions.
The choice of the sequence $\{x_n\}_{n\in\dN}$ in Section~\ref{sec:proof1} easily adjusts to possibly different choices of the sizes of the holes. It is not hard to see that in the construction in Section~\ref{sec:proof1} on the $n$th interconnection we are allowed to take any $\dl_n\in (0,\dl_n^\star)$ for some
sufficiently small $\dl_n^\star > 0$, which depends on the values of the parameters chosen on the previous steps of the construction. For the convenience of the reader we split  the remaining argument into two steps.
\subsubsection*{Step 1: modified sizes of the windows.}
Without loss of generality we can assume that $\dl_n \le a_1$ for all $n\in\dN$.
Since $\Omg_{n+1}\arr\Omg_n\cup\omg_{n+1}$ as $\dl_n\arr0$ in the sense of Definition~\ref{dfn:Mosco}, applying Proposition~\ref{prop:norm_resolvent} we conclude that for any $\eps >0$ there exists $\dl_n \in (0,\dl_n^\star)$ such that
\begin{equation}\label{eq:resdiffnorm}
\left\|(-\Delta_{\rm D}^{\Omg_{n+1}}+1)^{-1}-
(-\Delta_{\rm D}^{\Omg_{n}\cup\omg_{n+1}}+1)^{-1}\right\|\le \eps\,;
\end{equation}
note that both operators $-\Delta_{\rm D}^{\Omg_{n+1}}$ and $-\Delta_{\rm D}^{\Omg_n\cup\omg_{n+1}}$ act in the same Hilbert space $L^2(\Omg_{n+1})$.
Let the power $\ell\in\dN$ be fixed.
By purely algebraic manipulations the resolvent power difference of the Dirichlet Laplacians on $\Omg_{n+1}$ and $\Omg_n\cup\omg_{n+1}$ acting in the same Hilbert space $L^2(\Omg_{n+1})$ can be represented by
\begin{equation}\label{eq:resdiff}
\begin{aligned}
& (-\Delta_{\rm D}^{\Omg_{n+1}}+1)^{-\ell}-
(-\Delta_{\rm D}^{\Omg_{n}\cup\omg_{n+1}}+1)^{-\ell} =\\
&\quad\!=\!
\sum_{k=1}^\ell 
(-\Delta_{\rm D}^{\Omg_{n+1}}+1)^{-(\ell-k)}
\left[
(-\Delta_{\rm D}^{\Omg_{n+1}}+1)^{-1}\!-\!
(-\Delta_{\rm D}^{\Omg_n\cup\omg_{n+1}}+1)^{-1}
\right]
(-\Delta_{\rm D}^{\Omg_n\cup\omg_{n+1}}+1)^{-(k-1)}\,.
\end{aligned}
\end{equation}
We conclude from the Weyl spectral asymptotics for the Dirichlet Laplacian  on a bounded open set (see \eg~\cite[Thm.~12.14]{S12}) that for any $k\in\{1,2,\dots,\ell\}$ at least one of the two operators $(-\Delta_{\rm D}^{\Omg_{n+1}}+1)^{-(\ell-k)}$
and $(-\Delta_{\rm D}^{\Omg_{n}\cup\omg_{n+1}}+1)^{-(k-1)}$ belongs to the trace class provided that $\ell$ is chosen to be sufficiently large.
The trace norm of this trace class resolvent power can be estimated from above by a constant $C_n > 0$ independent of $k$ and $\dl_n$ while the operator norm of the other resolvent power can be estimated from above by $1$ thanks to non-negativity of the Dirichlet Laplacian. Let us comment on the possibility to get an upper bound on this trace norm independent of $\dl_n$. Indeed the operator $-\Delta^{\Omg_n\cup\omg_{n+1}}_{\rm D}$ is independent of $\dl_n$
and hence the trace norm of the corresponding resolvent power (provided it belongs to the trace class) is also independent of $\dl_n$. The trace norm of the resolvent power corresponding to $-\Delta^{\Omg_{n+1}}_{\rm D}$ (provided it belongs to the trace class) can be expressed as
\begin{equation}\label{eq:tracenorm2}
	\|(-\Delta^{\Omg_{n+1}}_{\rm D} +1)^{-(\ell-k)}\|_1 = \sum_{m=1}^\infty\frac{1}{
		\big[1+\lm_m(-\Delta^{\Omg_{n+1}}_{\rm D})\big]^{\ell-k}}
\end{equation}
where $\{\lm_m(-\Delta_{\rm D}^{\Omg_{n+1}})\}_{m\ge 1}$ are the eigenvalues $-\Delta_{\rm D}^{\Omg_{n+1}}$ enumerated in the non-decreasing way and repeated with multiplicities taken into account. It is easy to see by the min-max principle that 
the eigenvalues of $-\Delta_{\rm D}^{\Omg_{n+1}}$ are non-increasing functions of $\dl_n$. Hence, the trace norm in~\eqref{eq:tracenorm2}
is a non-decreasing function of $\dl_n$ and therefore it can be estimated from above by the respective trace norm for $\dl_n = a_1$.  

Combining the estimate~\eqref{eq:resdiffnorm} with the property of singular values in~\eqref{eq:sk1} and the definition of the trace norm in~\eqref{eq:tracenorm} we derive from~\eqref{eq:resdiff} the bound
\begin{equation}\label{eq:resdifftracenorm}
\left\|(-\Delta_{\rm D}^{\Omg_{n+1}}+1)^{-\ell}-
(-\Delta_{\rm D}^{\Omg_{n}\cup\omg_{n+1}}+1)^{-\ell}\right\|_1 \le C_n\ell \eps\,,\qquad n\in\dN\,. 
\end{equation}
It follows from~\eqref{eq:resdiffnorm} and~\eqref{eq:resdifftracenorm} that by choosing $\eps > 0$ sufficiently small we can find a sufficiently small $\dl_n \in (0,\dl_n^\star)$ so that
\begin{subequations}
\begin{equation}\label{eq:condnorm}
	\left\|(-\Delta_{\rm D}^{\Omg_{n+1}}+1)^{-1}-
	(-\Delta_{\rm D}^{\Omg_{n}\cup\omg_{n+1}}+1)^{-1}\right\| \le 2^{-n},\qquad n\in\dN\,,\\
\end{equation}	
\begin{equation}\label{eq:condtracenorm}
	\left\|(-\Delta_{\rm D}^{\Omg_{n+1}}+1)^{-\ell}-
	(-\Delta_{\rm D}^{\Omg_{n}\cup\omg_{n+1}}+1)^{-\ell}\right\|_1\le 2^{-n}\,,\qquad n\in\dN\,.
\end{equation}
\end{subequations} 
Thus, we have selected the sequence $\{\dl_n\}_{n\ge1}$ which fits the construction in Section~\ref{sec:proof1} and simultaneously the above two estimates hold.  

\subsubsection*{Step 2: trace class resolvent power difference}
Consider the sequence of quasi-conical unbounded domains
\[
\begin{aligned}
\Omg_1' & := \Omg,\\
\Omg_2' & := \Omg\cup Q_{\dl_1}(x_1+a_1),\\
\Omg_3' & := \Omg\cup Q_{\dl_1}(x_1+a_1)\cup Q_{\dl_2}(x_2+a_2),\\
\vdots &\\
\Omg_n' &:= \Omg\cup\Big(\bigcup_{k=1}^{n-1} Q_{\dl_k}(x_k+a_k)\Big),\\
\vdots &
\end{aligned}
\]
The Dirichlet Laplacian on $\Omg_n'$ can be decomposed into the orthogonal sum
\[
-\Delta_{\rm D}^{\Omg_n'} = (-\Delta_{\rm D}^{\Omg_n})\oplus (-\Delta_{\rm D}^{\Omg\sm\ov{\Omg_n}})
\]
with respect to the decomposition $L^2(\Omg) = L^2(\Omg_n)\oplus L^2(\Omg\sm\ov{\Omg_n})$.
In view of this decomposition we obtain the following representation for the resolvent difference of the Dirichlet Laplacians on~$\Omg_n'$ and~$\Omg_{n+1}'$:
\[
(-\Delta_{\rm D}^{\Omg_{n+1}'}+1)^{-1}
-
(-\Delta_{\rm D}^{\Omg_n'}+1)^{-1} = 
\left((-\Delta_{\rm D}^{\Omg_{n+1}}+1)^{-1}-
(-\Delta_{\rm D}^{\Omg_{n}\cup\omg_{n+1}}+1)^{-1}\right)\oplus 0,
\] 
where the orthogonal sum on the right-hand side is with respect to the decomposition $L^2(\Omg) = L^2(\Omg_{n+1}) \oplus L^2(\Omg\sm\ov{\Omg_{n+1}})$.
Hence, it follows from~\eqref{eq:condnorm} that the sequence $\dN\ni n\mapsto(-\Delta_{\rm D}^{\Omg_n'}+1)^{-1}$ is a Cauchy sequence in the operator norm. Thus, there exists an operator $\sfT$ in the Hilbert space $L^2(\Omg)$ such that
\begin{equation}\label{eq:limit}
	\left\|(-\Delta_{\rm D}^{\Omg_n'}+1)^{-1} -\sfT\right\| \arr 0\,,\qquad n\arr\infty\,.
\end{equation}
On the other hand the sequence of operators $\{-\Delta_{\rm D}^{\Omg_n'}\}_{n\ge 1}$ is decreasing in the form sense. In view of the inclusions $C^\infty_0(\Omg')\subset\cup_{n\ge1} H^1_0(\Omg_n')\subset H^1_0(\Omg')$
we get by \cite[Thm.~S.16]{RS1} that $-\Delta_{\rm D}^{\Omg_n'}$ converges to $-\Delta_{\rm D}^{\Omg'}$ in the strong resolvent sense. Thanks to~\eqref{eq:limit} we obtain that $\sfT = 
(-\Delta_{\rm D}^{\Omg'}+1)^{-1}$ and that 
$-\Delta_{\rm D}^{\Omg_n'}$ converges to $-\Delta_{\rm D}^{\Omg'}$ in the norm resolvent sense. 

 Let us introduce the sequence of the operators in the Hilbert space $L^2(\Omg')$ by
\[
	\sfS_n := (-\Delta_{\rm D}^\Omg+1)^{-\ell} - (-\Delta_{\rm D}^{\Omg_{n+1}'}+1)^{-\ell}\,,\qquad n\in\dN\,.
\]	
These operators can be represented as	
\[
\begin{aligned}
	\sfS_n & = \sum_{k=1}^n\Big[(-\Delta_{\rm D}^{\Omg_k'}+1)^{-\ell} - (-\Delta_{\rm D}^{\Omg_{k+1}'}+1)^{-\ell}\Big] \\
	&=\sum_{k=1}^n
	\Big[(-\Delta_{\rm D}^{\Omg_k\cup\omg_{k+1}}+1)^{-\ell} - (-\Delta_{\rm D}^{\Omg_{k+1}}+1)^{-\ell}\Big]\oplus 0\,.
\end{aligned}	
\]	
By~\eqref{eq:condtracenorm},
for all $n\in\dN$,
the operator $\sfS_n$ belongs to the trace class as a finite sum of trace class operators.
It follows from the norm resolvent convergence of $-\Delta_{\rm D}^{\Omg_n'}$ to $-\Delta_{\rm D}^{\Omg'}$ shown above  that
$\sfS_n$ converges to
$(-\Delta_{\rm D}^\Omg+1)^{-\ell} - (-\Delta_{\rm D}^{\Omg'}+1)^{-\ell}$ in the operator norm. Moreover, for any $n,m\in\dN$ with $m\ge n$ we find using the triangle inequality combined with~\eqref{eq:condtracenorm} that
\[
	\|\sfS_m-\sfS_n\|_1\le \sum_{k=n+1}^m	\left\|(-\Delta_{\rm D}^{\Omg_{k}\cup\omg_{k+1}}+1)^{-\ell}-(-\Delta_{\rm D}^{\Omg_{k+1}}+1)^{-\ell} \right\|_1 
	\le \sum_{k={n+1}}^m 2^{-k}\le 2^{-n}\,.
\]
Thus, $\{\sfS_n\}_{n\ge1}$ is a Cauchy sequence in the trace norm. Hence, it
converges in the trace norm to a trace class operator $\sfS$ acting in the Hilbert space $L^2(\Omg')$. Since the trace norm is stronger than the operator norm, we conclude that $\sfS_n$ converges to $\sfS$ also in the operator norm
and by the uniqueness of the limit we obtain that $\sfS = (-\Delta_{\rm D}^\Omg+1)^{-\ell} - (-\Delta_{\rm D}^{\Omg'}+1)^{-\ell}$.
Thus, the resolvent power difference
$(-\Delta_{\rm D}^{\Omg}+1)^{-\ell} - (-\Delta_{\rm D}^{\Omg'}+1)^{-\ell}$ is a trace class operator. Hence, by \cite[Thm.~X.4.8]{Kato} we have
\[
	\s_{\rm ac}(-\Delta_{\rm D}^{\Omg'}) = \s_{\rm ac}(-\Delta_{\rm D}^{\Omg}) =\varnothing.
\]
In the course of this construction we have also satisfied all the assumption in Section~\ref{sec:proof1}. Hence, the operator $-\Delta_{\rm D}^{\Omg'}$ possesses also a positive embedded eigenvalue $\lm_\infty$.

\section{Discussion}\label{sec:discussion}
It remains an open question whether it is possible to select the sizes of the windows so small that the singular continuous spectrum of the operator $-\Delta_{\rm D}^{\Omg'}$ is empty. If it
is possible to do, then we would get the Dirichlet Laplacian on a connected quasi-conical open set having only purely point spectrum, 
which densely fills the non-negative semi-axis. 

We also remark that our strategy can be adjusted to construct a similar connected quasi-conical open set having any finite number of (embedded) positive eigenvalues for the corresponding Dirichlet Laplacian.
This fact indicates that the point spectrum of the Dirichlet Laplacian on such domains is typically not limited by a single point.
\subsection*{Acknowledgement}
D.K. was supported by the EXPRO grant No. 20-17749X
of the Czech Science Foundation.
V.L. was supported by the grant No.~21-07129S 
of the Czech Science Foundation.

%
\bibliography{bib_final}

\providecommand{\bysame}{\leavevmode\hbox to3em{\hrulefill}\thinspace}
\providecommand{\MR}{\relax\ifhmode\unskip\space\fi MR }
\providecommand{\MRhref}[2]{%
  \href{http://www.ams.org/mathscinet-getitem?mr=#1}{#2}
}
\providecommand{\href}[2]{#2}
\begin{thebibliography}{10}

\bibitem{ABGM}
Y.~Avishai, D.~Bessis, B.~G. Giraud, and G.~Mantica, \emph{Quantum bound states
  in open geometries}, Phys.~Rev. B~ \textbf{44} (1991), 8028--8034.

\bibitem{Bonnet-Fliss-Hazard-Tonnoir_2016}
A.-S. Bonnet-Ben~Dhia, S.~Fliss, C.~Hazard, and A.~Tonnoir, \emph{A {R}ellich
  type theorem for the {H}elmholtz equation in a conical domain}, C. R. Acad.
  Sci. Paris \textbf{354} (2016), 27--32.

\bibitem{CK2}
L.~Cossetti and D.~Krej\v{c}i\v{r}\'ik, \emph{Absence of eigenvalues of
  non-self-adjoint {R}obin {L}aplacians on the half-space}, Proc. London. Math.
  Soc. \textbf{121} (2020), 584--616.

\bibitem{D'Ancona-Racke_2012}
P.~D'Ancona and R.~Racke, \emph{Evolution equations on non-flat waveguides},
  Arch. Ration. Mech. Anal. \textbf{206} (2012), 81--110.

\bibitem{Daners}
D.~Daners, \emph{Dirichlet problems on varying domains}, J.~Differential
  Equations \textbf{188} (2003), 591--624.

\bibitem{Daners_2008}
\bysame, \emph{Domain perturbation for linear and semi-linear boundary value
  problems}, Handbook of differential equations: Stationary partial
  differential equations. Vol. VI, Amsterdam: Elsevier/North Holland, 2008,
  pp.~1--81.

\bibitem{DLR12}
M.~Dauge, Y.~Lafranche, and N.~Raymond, \emph{Quantum waveguides with corners},
  ESAIM, Proc. \textbf{35} (2012), 14--45.

\bibitem{Edmunds-Evans}
D.~E. Edmunds and W.~D. Evans, \emph{Spectral theory and differential
  operators}, Oxford University Press, Oxford, 1987.

\bibitem{ESS1}
P.~Exner, P.~{\v S}eba, and P.~{\v S}{\v{t}}ov{\'\i}{\v c}ek, \emph{On
  existence of a bound state in an {L}-shaped waveguide}, Czech. J. Phys. B~
  \textbf{39} (1989), 1181--1191.

\bibitem{Fraenkel_1979}
L.~E. Fraenkel, \emph{On regularity of the boundary in the theory of {S}obolev
  spaces}, Proc. Lond. Math. Soc. \textbf{39} (1979), 385--427.

\bibitem{Glazman}
I.~M. Glazman, \emph{Direct methods of qualitative spectral analysis of
  singular differential operators}, Israel Program for Scientific Translations,
  1965.

\bibitem{Gohberg-Krein_1969}
I.~C. Gohberg and M.~G. Kre\v{\i}n, \emph{Introduction to the theory of linear
  nonselfadjoint operators in {H}ilbert space}, Amer. Math. Soc., 1969.

\bibitem{HSS91}
R.~Hempel, L.~A. Seco, and B.~Simon, \emph{{The essential spectrum of Neumann
  Laplacians on some bounded singular domains}}, {J. Funct. Anal.} \textbf{102}
  (1991), 448--483.

\bibitem{Jones_1953}
D.~S. Jones, \emph{The eigenvalues of {$\nabla^2 u + \lambda u = 0$} when the
  boundary conditions are given on semi-infinite domains}, Mathematical
  Proceedings of the Cambridge Philosophical Society \textbf{49} (1953), no.~4,
  668--684.

\bibitem{Kato}
T.~Kato, \emph{Perturbation theory for linear operators}, Springer-Verlag,
  Berlin, 1966.

\bibitem{RT}
J.~Rauch and M.~Taylor, \emph{Potential and scattering theory on wildly
  perturbed domains}, J.~Funct.~Anal. \textbf{18} (1975), 27--59.

\bibitem{RS1}
M.~Reed and B.~Simon, \emph{Methods of modern mathematical physics,
  {I}.~{F}unctional analysis}, Academic Press, New York, 1972.

\bibitem{Rellich_1943}
F.~Rellich, \emph{{\"U}ber das asymptotische {V}erhalten der {L}{\"o}sungen von
  {$\Delta u + \lambda u = 0$} in unendlichen {G}ebieten}, Jber. dtsch. MatVer.
  \textbf{53} (1943), 57--65.

\bibitem{S12}
K.~Schm{\"u}dgen, \emph{{Unbounded self-adjoint operators on Hilbert space}},
  vol. 265, Dordrecht: Springer, 2012.

\bibitem{Simon05}
B.~Simon, \emph{{Trace ideals and their applications}}, vol. 120, Providence,
  RI: American Mathematical Society, 2005.

\bibitem{Stollmann}
P.~Stollmann, \emph{A convergence theorem for {D}irichlet forms with
  applications to boundary value problems with varying domains}, Math.~Z.
  \textbf{219} (1995), 275--287.

\bibitem{Witsch_1990b}
K.~J. Witsch, \emph{Examples of embedded eigenvalues for the
  {D}irichlet-{L}aplacian in domains with infinite boundaries}, Math. Methods
  Appl. Sci. \textbf{12} (1990), 177--182.

\bibitem{Witsch_1990a}
\bysame, \emph{Examples of embedded eigenvalues for the {D}irichlet {L}aplacian
  in perturbed waveguides}, Math. Methods Appl. Sci. \textbf{12} (1990),
  91--93.

\end{thebibliography}
\bibliographystyle{amsplain}

\end{document}